\documentclass{ifacconf}

\usepackage{graphicx}      
\usepackage{epstopdf}
\usepackage{natbib}        

\usepackage{amsmath}
\usepackage{amssymb}
\usepackage{mathtools}
\mathtoolsset{showonlyrefs}

\usepackage{algorithm}
\usepackage{algpseudocode}

\begin{document}
\begin{frontmatter}

\title{Successive Convexification of Non-Convex Optimal Control Problems with 
	State Constraints} 


\author[First]{Yuanqi Mao} 
\author[First]{Daniel Dueri}
\author[First]{Michael Szmuk}
\author[First]{Beh\c{c}et A\c{c}\i kme\c{s}e}

\address[First]{University of Washington, Seattle, WA 98105 USA 
	(e-mail: {\tt yqmao@uw.edu}, {\tt dandueri@uw.edu}, {\tt mszmuk@uw.edu}, {\tt behcet@uw.edu}).}

\begin{abstract}                
This paper presents a Successive Convexification (\texttt{SCvx}) algorithm to solve a class of non-convex optimal control problems with certain types of state constraints. Sources of non-convexity may include nonlinear dynamics and non-convex state/control constraints. To tackle the challenge posed by non-convexity, first we utilize exact penalty function to handle the nonlinear dynamics. Then the proposed algorithm successively convexifies the problem via a \textit{project-and-linearize} procedure. Thus a finite dimensional convex programming subproblem is solved at each succession, which can be done efficiently with fast Interior Point Method (IPM) solvers. Global convergence to a local optimum is demonstrated with certain convexity assumptions, which are satisfied in a broad range of optimal control problems. The proposed algorithm is particularly suitable for solving trajectory planning problems with collision avoidance constraints. Through numerical simulations, we demonstrate that the algorithm converges reliably after only a few successions. Thus with powerful IPM based custom solvers, the algorithm can be implemented onboard for real-time autonomous control applications.
\end{abstract}

\begin{keyword}
Optimal control, State constraints, Convex optimization, Successive convexification, Autonomous systems, Trajectory optimization
\end{keyword}

\end{frontmatter}

\section{Introduction}

Non-convex optimal control problems emerge in a broad range of science and engineering disciplines. Finding a global solution to these problems is generally considered NP-hard. Heuristics like simulated annealing, \cite{bertsimas2010robust}, or combinatorial methods like mixed integer programming, \cite{richards2002}, can compute globally optimal solutions for special classes of problems. In many engineering applications however, finding a local optimum or even a feasible solution with much less computational effort is a more favorable route. This is particularly the case with real-time control systems, where efficiency and convergence guarantees are more valuable than optimality. An example in aerospace applications is the planetary landing problem, see~\cite{pointing2013, lars_sys12, steinfeldt2010guidance}. Non-convexities in this problem include minimum thrust constraints, nonlinear gravity fields and nonlinear aerodynamic forces. State constraints can also render the problem non-convex. A classic example is imposing collision avoidance constraints. For instance, \cite{accikmese2006convex} discusses the collision avoidance in formation reconfiguration of spacecraft, \cite{augugliaro2012generation} considers the generation of collision-free trajectories for a quad-rotor fleet, and~\cite{liu2014solving} study the highly constrained rendezvous problem.

Given the complexity of such non-convex problems, traditional Pontryagin's maximum principle-based approaches, e.g.~\cite{rockafellar1972state}, can fall short. On the other hand, directly applying optimization methods to solve the discretized optimal control problems has gained in popularity thanks to algorithmic advancements in nonlinear programming, see e.g.~\cite{hull1997, buskens}. However, general nonlinear optimization can sometimes be intractable in the sense that a bad initial guess could potentially lead to divergence, and also there are few known bounds on the computational effort needed to reach optimality. This makes it difficult to implement for real-time or mission critical applications because these applications cannot afford either divergence or a large amount of computational effort. Convex optimization, on the other hand, can be reliably solved in polynomial time to global optimality, see e.g.~\cite{BoydConvex}. More importantly, recent advances have shown that these problems can be solved in real-time by both generic Second Order Cone Programming (SOCP) solvers, e.g.~\cite{alexd}, and by customized solvers which take advantage of specific problem structures, e.g.~\cite{mattingley2012, dueri2014automated}. This motivates researchers to formulate optimal control problems in a convex programming framework for real-time purposes, e.g., real-time Model Predictive Control (MPC), see~\cite{houska2011auto, Zeilinger2014683}.

In order to take advantage of these powerful convex programming solvers, one crucial step is to convexify the originally non-convex problems. Recent results on a technique known as \textit{lossless convexification}, e.g.~\cite{behcet_aut11, lars_sys12, matt_aut1} have proven that certain types of non-convex control constraints can be posed as convex ones without introducing conservativeness. \cite{liu2015entry} also gives a result on convexification of control constraints for the entry guidance application. For nonlinear dynamics and non-convex state constraints, collision avoidance constraints in particular, one simple solution is to perform query-based collision checks, see \cite{allen2016real}. However, to be more mathematically tractable, \cite{augugliaro2012generation, schulman2014motion,  chen2015decoupled} propose to use (variations of) sequential convex programming (SCP) to iteratively convexify non-convexities. While these methods usually perform well in practice, no convergence results have been reported yet. As an effort to tackle this open challenge, \cite{SCvx_cdc16} propose an successive convexification (\texttt{SCvx}) algorithm that successively convexifies the dynamics with techniques like \textit{virtual control} and \textit{trust regions}, and more importantly, give a convergence proof of that algorithm.

To include state constraints in the \texttt{SCvx} algorithmic framework, a few enhancements need to be made. \cite{hauser2006barrier} relax the state constraints by using a barrier function, but do not provide theoretical guarantees. In this paper, we relax the dynamic equations into inequalities by using an exact penalty function, and then we propose a \textit{project-and-linearize} procedure to handle both the state constraints and the relaxed dynamics. While introducing conservativeness is inevitable in the process, the proposed algorithm preserves much more feasibility than similar results in~\cite{rosen1966iterative, liu2014solving}. Finally under some mild assumptions, we present a convergence proof, which not only guarantees that the algorithm will converge, but also demonstrates that the convergent solution recovers local optimality for the original problem.

One clear advantage of the algorithm proposed in this paper is that this algorithm does not have to resort to trust regions, as in e.g.~\cite{SCvx_cdc16,szmuk2016successive,szmuk2017successive}, to guarantee convergence. This property allows the algorithm to potentially take a large step in each succession, thereby greatly accelerating the convergence process, which is exactly the case shown by the numerical simulations. It also worth noting that the proposed algorithm only uses the Jacobian matrix, i.e. first-order information; therefore, we do not have to compute Hessian matrices, otherwise that task itself could be computationally expensive. To the best of our knowledge, the main contributions of this work are:
\begin{itemize}
	\item An extended \texttt{SCvx} algorithm with \textit{project-and-linearize} to handle both nonlinear dynamics and non-convex state constraints.
	\item A convergence proof with local optimality recovery.
\end{itemize}


\section{Successive Convexification}


\subsection{Problem Formulation}

In this paper, we consider the following discrete optimal control problem:
\begin{subequations}
	\begin{equation}
	\min J(x_{i},u_{i}):= \sum_{i=1}^{T}\phi(x_{i},u_{i}), \label{cost:original}
	\end{equation}
	subject to
	\begin{align}
	&x_{i+1}-x_{i}=f(x_{i},u_{i}) &i=1,2,T-1, \label{eq:dynamics} \\
	&h(x_{i})\geq 0 &i=1,2,T, \label{eq:state_con} \\
	&u_{i}\in U_{i} \subseteq \mathbb{R}^{m} &i=1,2,T-1, \\
	&x_{i}\in X_{i} \subseteq \mathbb{R}^{n} &i=1,2,T.
	\end{align}
\end{subequations}
Here, $ x_{i},u_{i} $ represent discrete state/control at each temporal point, $ T $ denotes the final time, and $ X_{i},U_{i} $ are assumed to be convex and compact sets. We also assume that the objective function in~\eqref{cost:original} is continuous and convex, as is the case in many optimal control applications. For example, the minimum fuel problem has $ \phi(x_{i},u_{i}) = \|u_{i}\| $, and the minimum time problem has $ \phi(x_{i},u_{i}) = 1 $. 
Equation~\eqref{eq:dynamics} represent the system dynamics, where $ f(x_{i},u_{i})\in \mathbb{R}^{n} $ is, in general, a nonlinear function that is at least twice differentiable.
Equation~\eqref{eq:state_con} are the additional state constraints, where $ h(x_{i})\in \mathbb{R}^{s} $ is also at least twice differentiable, and could be nonlinear as well. Note that we do not impose non-convex control constraints here because we can leverage lossless convexification, see e.g.~\cite{behcet_aut11}, to convexify them beforehand. Finally, note that ~\eqref{eq:dynamics} and~\eqref{eq:state_con} render the problem non-convex.

\begin{figure}[!ht]
	\begin{center}
		\includegraphics[width=0.25\textwidth]{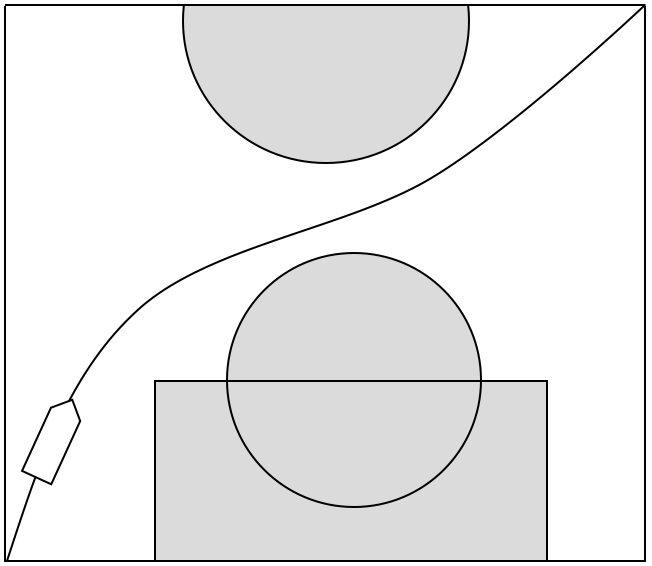}
		\caption{Convex shaped keep-out zone as state constraints.} 
		\label{fig:keep_out}
	\end{center}
\end{figure}

To facilitate the convergence proof, we need a few more assumptions on $ f(x_{i},u_{i}) $ and $ h(x_{i}) $. First, we assume each component of $ f(x_{i},u_{i}) $ is a convex function over $ x_{i} $ and $ u_{i} $. In fact, a wide range of optimal control applications, for example systems with double integrator dynamics and aerodynamic drag (constant speed), satisfy this assumption. It also includes all linear systems. Similarly, we assume each component of $ h(x_{i}) $ is convex. An example for this assumption is collision avoidance constraints, where the shape of each keep-out zone is convex or can be convexly decomposed. See Figure~\ref{fig:keep_out} for a simple illustration. Note that at this stage, these convexity assumptions do not change the non-convex nature of the problem.

To formalize, we give the following hypotheses:

\begin{hypo} \label{hypo:convex_f}
	$ f_{j}(x_{i},u_{i}) $ is a convex function over $ x_{i} $ and $ u_{i} $, $ \forall j = 1,2,...n $.
\end{hypo}

\begin{hypo} \label{hypo:convex_h}
	$ h_{j}(x_{i}) $ is a convex function over $ x_{i} $, $ \forall j = 1,2,...s $.
\end{hypo}

To convert the optimal control problem into a finite dimensional optimization problem, we treat state variables $ x_{i} $ and control variables $ u_{i} $ as a single variable $ y = \left( x_{1}^{T},..., x_{T}^{T}, u_{1}^{T},..., u_{T-1}^{T} \right)^{T} \in \mathbb{R}^{N} $, where $ N = m(T-1)+nT $. Let $ Y $ be the Cartesian product of all $ X_{i} $ and $ U_{i} $, then $ y \in Y $. $ Y $ is a convex and compact set because $ X_{i} $ and $ U_{i} $ are. In addition, we let $ g_{i}(x_{i}, u_{i}) = f(x_{i}, u_{i}) - x_{i+1} + x_{i} $, and $ g(y) = \left( g_{1}^{T}, g_{2}^{T},..., g_{T-1}^{T} \right)^{T} \in \mathbb{R}^{n(T-1)} $. Note that each component of $ g(y) $ is convex over $ y $ by Hypothesis~\ref{hypo:convex_f}, and the dynamic equation~\eqref{eq:dynamics} becomes $ g(y) = 0 $. We also let $ h(y) = \left( h(x_{1})^{T}, h(x_{2})^{T},..., h(x_{T})^{T} \right)^{T} \in \mathbb{R}^{sT} $, then each component of $ h(y) $ is convex over $ y $ by Hypothesis~\ref{hypo:convex_h}, and~\eqref{eq:state_con} becomes $ h(y)\geq 0 $. In summary, we have the following non-convex optimization problem:
\begin{equation} \label{eq:ori_prob}
J(y^{*}) = \underset{y}{min}\{J(y)\ |\ y\in Y, g(y)=0, h(y)\geq 0 \}.
\end{equation}
By leveraging the theory of exact penalty methods, see e.g.~\cite{han1979exact},~\cite{SCvx_cdc16}, we move $ g(y)=0 $ into the objective function without compromising optimality:

\begin{thm}[\textbf{Exactness}] \label{thm:exactness}
	Let $ P(y) = J(y)+\lambda \|g(y)\|_{1} $ be the penalty function, and $ \bar{y} $ be a stationary point of
	\begin{equation} \label{eq:pen_prob}
	\underset{y}{min}\{P(y)\ |\ y\in Y, g(y)\geq 0, h(y)\geq 0 \}
	\end{equation}
	with Lagrangian multiplier $ \bar{\mu} $ for equality constraints. Then, if the penalty weight $ \lambda $ satisfies $ \lambda\geq \|\bar{\mu}\|_{\infty} $, and if $ \bar{y} $ is feasible for~\eqref{eq:ori_prob}, then $ \bar{y} $ is a critical point of~\eqref{eq:ori_prob}.
\end{thm}
Since each component of $ g(y) $ is convex, and $ \|\cdot \|_{1} $ is convex and nondecreasing due to the constraint $ g(y)\geq 0 $, then $ P(y) $ is a convex function by the composition rule of convex functions, see~\cite{BoydConvex}. This marks our first effort towards the convexification of~\eqref{eq:ori_prob}.

Let $ q(y) = \left( g(y)^{T}, h(y)^{T}\right)^{T} \in \mathbb{R}^{M} $, where $ M = sT+n(T-1) $, then we may rewrite~\eqref{eq:pen_prob} as
\begin{equation} \label{eq:final_prob}
\underset{y}{min}\{P(y)\ |\ y\in Y, q(y)\geq 0 \}
\end{equation}
where $ P(y) $ is continuous and convex, and each component of $ q(y) $ is a convex function. By doing this, we are essentially treating constraints due to dynamics as another keep-out zone. Denote
\begin{equation*}
F = \{y\ |\ y\in Y, q(y)\geq 0 \}
\end{equation*}
as the feasible set. Note that $ F $ is compact but not convex.

\subsection{Project and Linearize}
Since both $ f(x_{i},u_{i}) $ and $ h(x_{i}) $ are at least twice differentiable, we have $ q(y) \in C^{2}(Y) $ as well. For any point $ z \in F $, let $ \triangledown_{y} q(z) $ be the Jacobian matrix of $ q(y) $ evaluated at $ z $. 

\begin{figure}[!ht]
	\begin{center}
		\includegraphics[width=0.3\textwidth]{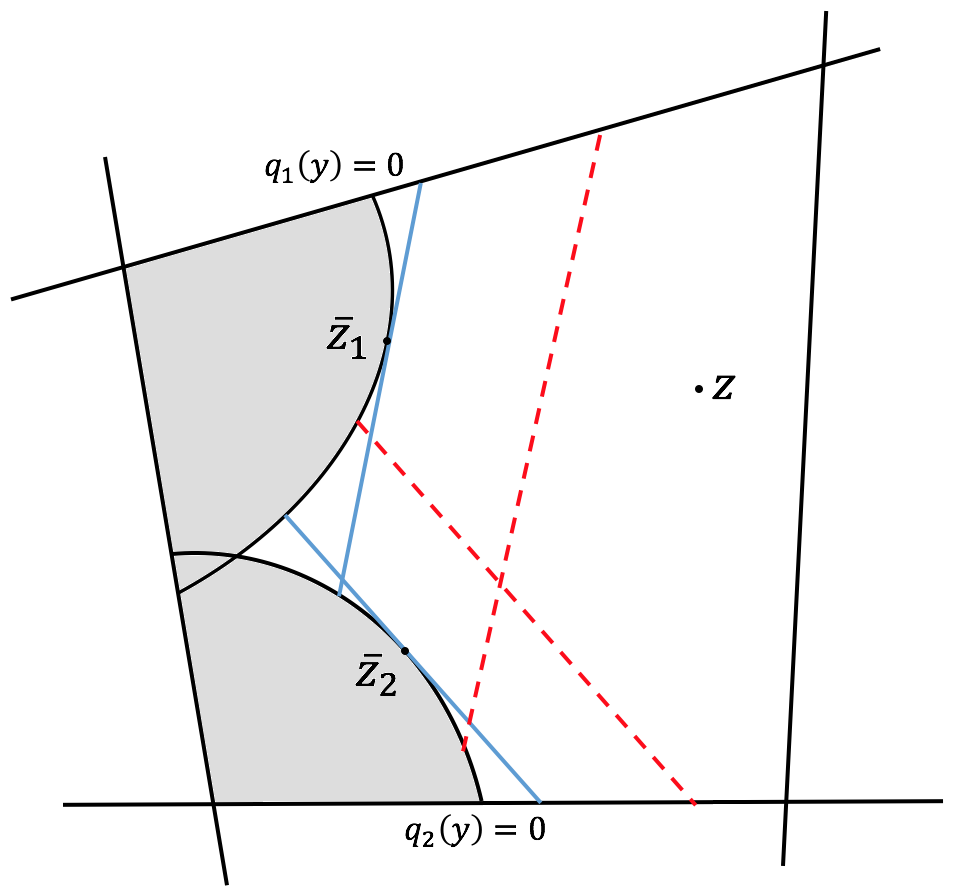}
		\caption{Project-and-linearize vs directly linearize.} 
		\label{fig:project}
	\end{center}
\end{figure}

Now, if we directly linearize $ q(y) $ at $ z $ as in~\cite{liu2014solving}, there may be a gap between the linearized feasible region and $ F $ since $ z $ could be in the interior of $ F $. The gap will increase as $ z $ moves further away from the boundary $ q(y)=0 $. This is not a desirable situation, because a fairly large area of the feasible region is not utilized. In other words, we introduced artificial conservativeness. To address this issue, we introduce a \textit{projection} step, which essentially projects $ z $ onto each constraint, and obtains each projection point. Then, we linearize each constraint at its own projection point. See Figure~\ref{fig:project} for an illustration in $ \mathbb{R}^{2} $.

To formalize, let $ q_{j}(y), j=1,2...,M $ represent each component of $ q(y) $, i.e. each constraint. Note that $ q_{j}(y) $ is a convex function, hence $ q_{j}(y) \leq 0 $ is a closed convex set. Using the well-known Hilbert projection theorem, see e.g.\cite{wulbert_project}, we have

\begin{thm}[\textbf{Uniqueness of projection}]
	For any $ z \in F $, there exists a unique point
	\begin{equation} \label{eq:project}
	\bar{z}_j = \underset{y}{argmin} \left\lbrace \|z-y\|_{2}\ |\ q_{j}(y) \leq 0 \right\rbrace,
	\end{equation}
	called the projection of $ z $ onto $ q_{j}(y) \leq 0 $.
\end{thm}
Equation~\eqref{eq:project} is a simple convex program of low dimension that can be solved quickly (sub-milliseconds) using any convex programming solver. Alternatively, for some special convex sets (e.g. cylinders), \eqref{eq:project} can be solved analytically, which is even faster. Doing this for each constraint, we obtain a set of projection points, $ \{\bar{z}_{1}, \bar{z}_{2},..., \bar{z}_{M} \} $. Note that these projection points must lie on the boundary of $ q_{j}(y) \leq 0 $, i.e.
\begin{equation*}
	q_{j}(\bar{z}_{j}) = 0, \quad \forall \; j = 1,..., M.
\end{equation*}
For a fixed $ z \in F $, let $ l_{j}(y,z) $ be the linear approximation of $ q_{j}(y) $:
\begin{equation} \label{eq:linearization}
l_{j}(y,z) = \triangledown_{y}q_{j}(\bar{z}_{j})(y - \bar{z}_{j}),
\end{equation}
and let $ l(y,z) = (l_{1}^{T}, l_{2}^{T},..., l_{M}^{T})^{T} \in \mathbb{R}^{M} $. For each $ z \in F $, denote
\begin{equation*}
	F_{z} = \{y\ |\ y\in Y,\, l(y,z)\geq 0 \}
\end{equation*}
as the feasible region after linearization. $ F_{z} $ also defines a point-to-set mapping, $ F_{z}: z \rightarrow F_{z} $.
Note that each component of $ l(y,z)\geq 0 $, i.e. $ l_{j}(y,z)\geq 0 $ represents a half-space. Hence $ l(y,z)\geq 0 $ is the intersection of half-spaces, which means $ F_{z} $ is a convex and compact set.

\subsubsection{Remark.} Convexification of $ F $ by using $ F_{z} $ also inevitably introduces conservativeness, but one can verify that it is the best we can do to maximize feasibility while preserving convexity.

The following lemma gives an invariance result regarding the point-to-set mapping $ F_{z} $. It is essential to our subsequent analyses.

\begin{lem}[\textbf{Invariance of $ F_{z} $}] \label{lem:invariance}
	For each $ z \in F $,
	\begin{equation*}
	z \in F_{z} \subseteq F.
	\end{equation*}
\end{lem}

\begin{pf}
	For each $ z \in F $, and $ \forall j = 1,2,..., M $, from~\eqref{eq:linearization}, we have
	\begin{equation*}
	l_{j}(z,z) = \triangledown_{y}q_{j}(\bar{z}_{j})(z - \bar{z}_{j}).
	\end{equation*}
	Since $ \bar{z}_{j} $ is the projection, $ (z - \bar{z}_{j}) $ is the normal vector at $ \bar{z}_{j} $, which is aligned with the gradient $ \triangledown_{y}q_{j}(\bar{z}_{j}) $. Hence $ \triangledown_{y}q_{j}(\bar{z}_{j})(z - \bar{z}_{j}) \geq 0 $, i.e. $ l_{j}(z,z) \ge 0, \forall j = 1,2,..., M $, i.e. $ z \in F_{z} $.
	
	Furthermore, since $ q_{j}(y) $ is a convex function, we have for any $ y \in F_{z} $,
	\begin{equation*}
	q_{j}(y) \geq q_{j}(\bar{z}_{j}) + \triangledown_{y}q_{j}(\bar{z}_{j})(y - \bar{z}_{j}) 
	= l_{j}(y,z) \geq 0,
	\end{equation*}
	which means $ y \in F $. Hence $ F_{z} \subseteq F $. \qed
\end{pf}

\subsection{The \texttt{SCvx} Algorithm}
Now that we have a convex and compact feasible region $ F_{z} $ and a convex objective function $ P(y) $, we are ready to present a successive procedure to solve the non-convex problem in~\eqref{eq:final_prob}. Note that the feasible region $ F_{z} $ is defined by $ z $. Therefore, if we start from a point $ z^{(0)} \in F $, a sequence $ \{z^{(k)} \} $ will be generated, where
\begin{equation} \label{eq:sub_prob}
z^{(k+1)} = \underset{y}{argmin}\{P(y)\ |\ y\in F_{z^{(k)}} \}, \quad k=0,1,...
\end{equation}
This is a convex programming sub-problem, whose global minimizer is attained at $ z^{(k+1)} $. At these intermediate steps, $ z^{(k+1)} $ may not be the optimal solution to~\eqref{eq:final_prob}. Our goal, however, is to prove that this sequence $ \{z^{(k)} \} $ converges to a limit point $ z^{*} $, and that this limit point solves~\eqref{eq:final_prob} by \textit{project-and-linearize} at $ z^{*} $ itself, i.e., it is a ``fixed-point" satisfying
\begin{equation} \label{eq:optimal}
z^{*} = \underset{y}{argmin}\{P(y)\ |\ y\in F_{z^{*}} \}.
\end{equation}
More importantly, we want to show that $ z^{*} $ gives a local optimum to~\eqref{eq:final_prob} convexified at $ z^{*} $ itself. Then by solving a sequence of convex programming subproblems, we effectively solved the non-convex optimal control problem in~\eqref{eq:ori_prob} because of Fact~\ref{thm:exactness}. Therefore we call this procedure the Successive Convexification (\texttt{SCvx}) algorithm. It is summarized in Algorithm~\ref{algo:SCvx}.
\begin{algorithm}
	\caption{Successive Convexification Algorithm}
	\label{algo:SCvx}
	\begin{algorithmic}[1]
		\Procedure{\texttt{SCvx}}{$z^{(0)},\lambda$}
		\State \textbf{input} Initial point $ z^{(0)} \in F $. Penalty weight $ \lambda \geq 0 $.
		\While {not converged}
		\State \textbf{step 1} At each succession $ k $, we have the current point $ z^{(k)} $. For each constraint $ q_{j}(y) $, solve the problem in~\eqref{eq:project} to get a projection point $ \bar{z_{j}^{(k)}} $.
		\State \textbf{step 2} Construct the convex feasible region $ F_{z^{(k)}} $ by using~\eqref{eq:linearization}.
		\State \textbf{step 3} Solve the convex subproblem in~\eqref{eq:sub_prob} to get $ z^{(k+1)} $. Let $ z^{(k)}\leftarrow z^{(k+1)} $ and go to the next succession.
		\EndWhile
		\State \textbf{return} $ z^{(k+1)} $.
		\EndProcedure
	\end{algorithmic}
\end{algorithm}

\section{Convergence Analysis}

In this section, we proceed to show that Algorithm~\ref{algo:SCvx} does converge to a point $ z^{*} $ that indeed satisfies~\eqref{eq:optimal}. First we must assume the application of regular constraint qualifications, namely the Linear Independence Constraint Qualification (LICQ) and the Slater's condition. They can be formalized as the following:

\begin{hypo}[\textbf{LICQ}] \label{hypo:licq}
	For each $ z \in F $, the Jacobian matrix $ \triangledown_{y} q(z) $ has full row rank, i.e. rank($ \triangledown_{y} q(z) $) $ = M $.
\end{hypo}

\begin{hypo}[\textbf{Slater's condition}] \label{hypo:slater}
	For each $ z \in F $, the convex feasible region $ F_{z} $ contains interior points.
\end{hypo}

\begin{figure}[!ht]
	\begin{center}
		\includegraphics[width=0.4\textwidth]{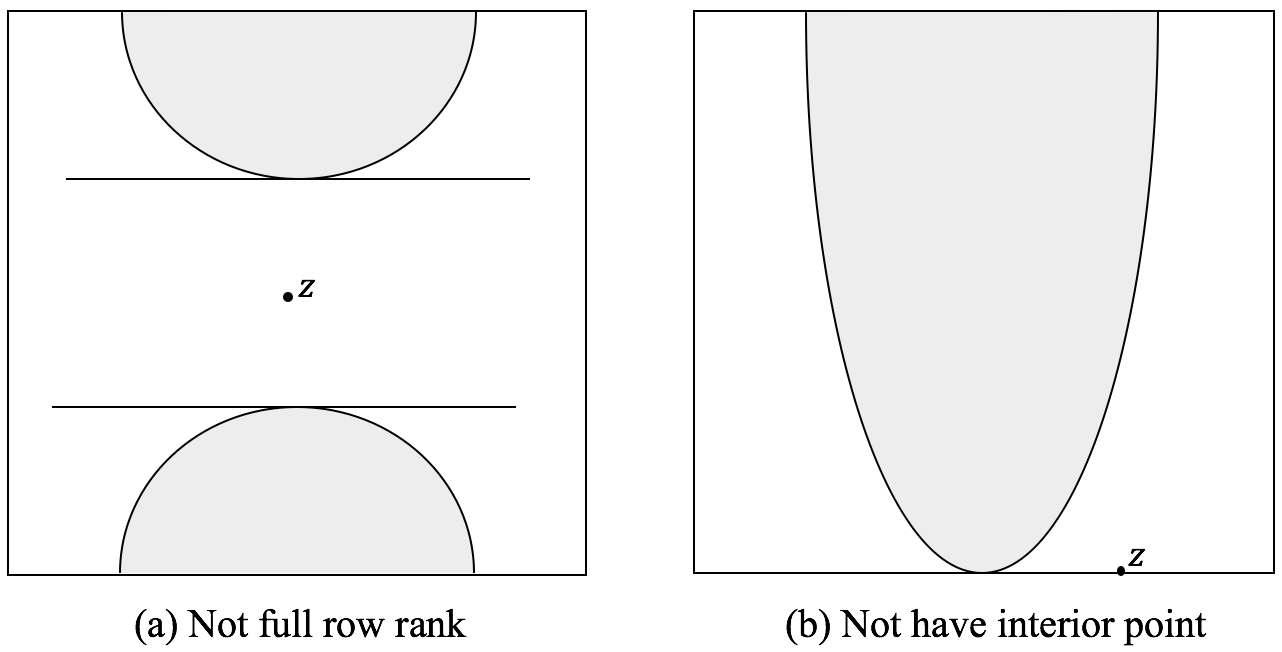}
		\caption{Cases where LICQ or Slater's condition fails} 
		\label{fig:licq}
	\end{center}
\end{figure}

These assumptions do impose some practical restrictions on our feasible region. Figure~\ref{fig:licq} shows some examples where LICQ or Slater's condition might fail. For scenarios like (a), we may perturb our discrete points to break symmetry, For scenarios like (b), we have to assume the connectivity of the feasible region. In other words, our feasible region cannot be degenerate at some point, for example, in collision avoidance the feasible region is not allowed to be completely obstructed.

In order to analyze the convergence, first we need to show that the point-to-set mapping $ F_{z} $ is continuous in the sense that given any point $ z^{(1)} \in F $ and $ y^{(1)} \in F_{z^{(1)}} $, then for any point $ z^{(2)} \in F $ in the neighborhood of $ z^{(1)} $, there exists a $ y^{(2)} \in F_{z^{(2)}} $ that is close enough to $ y^{(1)} $.

First, We have the following lemma:
\begin{lem}[\textbf{Lipschitz continuity of $ l(y,z) $}] \label{lem:lipschitz}
	The linear approximation $ l(y,z) $ in~\eqref{eq:linearization} is Lipschitz continuous in $ z $, that is
	\begin{equation} \label{eq:lipschitz}
	\| l(y,z^{(1)}) - l(y,z^{(2)}) \| \leq \gamma \| z^{(1)} - z^{(2)} \|,
	\end{equation}
	 for any $ z^{(1)}, z^{(2)} \in F $ and constant $ \gamma $.
\end{lem}

\begin{pf}
	Each $ l_{j}(y,z) $ is in fact a composition of two mappings. The first mapping maps $ z $ to its projection $ \bar{z}_{j} $. This mapping is defined by the optimization problem in~\eqref{eq:project}. It is a well-known result that this mapping is non-expansive (i.e. Lipschitz continuous with constant 1), See e.g.~\cite{wulbert_project}. The second mapping is defined by the auxiliary function:
	\begin{equation*}
	a_{j}(y,\bar{z}_{j}) = \triangledown_{y}q_{j}(\bar{z}_{j})(y - \bar{z}_{j}),
	\end{equation*}
	Since $ q(y) \in C^{2}(Y) $, $ a_{j}(y,\bar{z}_{j}) $ is Lipschitz continuous in $ \bar{z}_{j} $.
	By the composition rule of two Lipschitz continuous functions, we have $ l_{j}(y,z) $ is Lipschitz continuous, for all $ j = 1,2,..., M $. Therefore, sum over $ j $ gives the Lipschitz continuity of $ l(y,z) $. \qed
\end{pf}

Now with Hypothesis~\ref{hypo:licq},~\ref{hypo:slater} and Lemma~\ref{lem:lipschitz}, we are ready to prove the continuity of point-to-set mapping $ F_{z} $. The result is given as follows:

\begin{lem}[\textbf{Continuity of mapping $ F_{z} $}] \label{lem:continuity}
	Given $ z^{(1)} \in F $ and $ y^{(1)} \in F_{z^{(1)}} $, then given $ \epsilon > 0 $, there exists a $ \delta > 0 $ so that for any point $ z^{(2)} \in F $ with $ \|z^{(2)} - z^{(1)}\| \leq \delta $, there exists a $ y^{(2)} \in F_{z^{(2)}} $ such that $ \|y^{(2)} - y^{(1)}\| \leq \epsilon $.
\end{lem}

\begin{pf}
	From Hypothesis~\ref{hypo:licq}, we know that the Jacobian matrix $ \triangledown_{y} q(z) $ has full row rank. Thus matrix $ \triangledown_{y} q(z)\triangledown_{y} q(z)^{T} $ is symmetric and positive definite for any $ z \in F $. Consequently, there exists $ \beta $ such that
	\begin{equation} \label{eq:deriv_ineq}
	\| (\triangledown_{y} q(z)\triangledown_{y} q(z)^{T})^{-1} \| \leq \beta^{2}, \quad \forall z \in F.
	\end{equation}
	If $ y^{(1)} \in F_{z^{(2)}} $, then take $ y^{(2)} = y^{(1)} $ such that $ \|y^{(2)} - y^{(1)}\| =0 \leq \epsilon $. Now suppose $ y^{(1)} \notin F_{z^{(2)}} $, then there exists at least one $ j $, such that $ l_{j}(y^{(1)}, z^{(2)}) < 0 $. Let
	\begin{equation*}
	\bar{l}_{j} = \begin{cases}
	l_{j}(y^{(1)}, z^{(2)}) & \quad l_{j}(y^{(1)}, z^{(2)}) < 0, \\
	0 & \quad l_{j}(y^{(1)}, z^{(2)}) \geq 0.
	\end{cases}
	\end{equation*}
	Note that $ l_{j}(y^{(1)}, z^{(1)}) \geq 0 $, then by definition,
	\begin{equation*}
	| \bar{l}_{j} | \leq | l_{j}(y^{(1)}, z^{(1)}) - l_{j}(y^{(1)}, z^{(2)}) |, \quad \forall j = 1, 2, ..., M.
	\end{equation*}
	Hence we have
	\begin{equation} \label{eq:norm_L_bar}
	\| \bar{l} \| \leq \| l(y^{(1)}, z^{(1)}) - l(y^{(1)}, z^{(2)}) \| \leq \gamma \| z^{(1)} - z^{(2)} \|.
	\end{equation}
	The second inequality follows from~\eqref{eq:lipschitz} in Lemma~\ref{lem:lipschitz}.
	
	Now we consider two cases: \\
	Case 1: $ y^{(1)} $ is an interior point of $ Y $, then $ \exists \epsilon_{1} \in (0, \epsilon] $, such that $ y \in Y $ for all $ y $ satisfies $ \| y - y^{(1)} \| \leq \epsilon_{1} $. \\
	Let $ \|z^{(2)} - z^{(1)}\| \leq \delta $, with $ \delta = \epsilon_{1}/\beta\gamma $, and let
	\begin{equation} \label{eq:y2}
	y^{(2)} = y^{(1)} - \triangledown_{y} q(z^{(2)})^{T} (\triangledown_{y} q(z^{(2)})\triangledown_{y} q(z^{(2)})^{T})^{-1} \bar{l}.
	\end{equation}
	We need to verify that $ y^{(2)} \in F_{z^{(2)}}$. First, we have
	\begin{equation} \label{eq:Lj2}
	l_{j}(y^{(2)},z^{(2)}) = \triangledown_{y}q_{j}(\bar{z}_{j}^{(2)})(y^{(2)} - \bar{z}_{j}^{(2)}).
	\end{equation}
	Substitute $ y^{(2)} $ from~\eqref{eq:y2} into the stacked form of~\eqref{eq:Lj2}, and then unstack, we get
	\begin{align*}
	l_{j}(y^{(2)},z^{(2)}) &= \triangledown_{y}q_{j}(\bar{z}_{j}^{(2)})(y^{(1)} - \bar{z}_{j}^{(2)}) - \bar{l}_{j} \\
	&= l_{j}(y^{(1)},z^{(2)}) - \bar{l}_{j} \\
	&\geq 0.
	\end{align*}
	The last inequality follows from the definition of $ \bar{l}_{j} $. Therefore, $ l(y^{(2)},z^{(2)}) \geq 0 $. \\
	Next we need to show that $ y^{(2)} \in Y $. Rearrange terms in~\eqref{eq:y2}, we have
	\begin{align*}
	\| y^{(2)} - y^{(1)} \|^{2} &= \bar{l}^{T} (\triangledown_{y} q(z^{(2)})\triangledown_{y} q(z^{(2)})^{T})^{-1} \bar{l} \\
	&\leq \beta^{2} \| \bar{l} \|^{2} \\
	&\leq \beta^{2} \gamma^{2} \| z^{(1)} - z^{(2)} \|^{2}
	\end{align*}
	The inequalities follow from~\eqref{eq:deriv_ineq} and~\eqref{eq:norm_L_bar}. Thus we have $ \| y^{(2)} - y^{(1)} \| \leq \beta\gamma\delta = \epsilon_{1} $, i.e. $ y^{(2)} \in Y $. So now we have verified $ y^{(2)} \in F_{z^{(2)}}$. Since $ \epsilon_{1} \leq \epsilon $, we have $ \| y^{(2)} - y^{(1)} \| \leq \epsilon $, so that $ F_{z} $ is continuous.
	
	Case 2: $ y^{(1)} $ is a boundary point of $ Y $. From Hypothesis~\ref{hypo:slater}, $ F_{z^{(1)}} $ has interior points. Also since $ F_{z^{(1)}} $ is a convex set, it is a well known fact that there are interior points in the $ \epsilon- $Neighborhood of every point in $ F_{z^{(1)}} $. Then we can apply the same argument as in Case 1, to get that $ F_{z} $ is continuous. \qed
\end{pf}

One way to show convergence is to demonstrate the convergence of objective functions, $ P(z^{(k)}) $. Now let's define
\begin{equation}
\Phi (z) \coloneqq \underset{y}{min}\{P(y)\ |\ y\in F_{z} \}
\end{equation}
to be the function that maps the point $ z $ we are solving at in each iteration to the optimal value of the objective function in $ F_{z} $. By using the continuity of the point-to-set mapping $ F_{z} $, the following lemma gives the continuity of the function $ \Phi (z) $.

\begin{lem} \label{lem:continuity_phi}
	$ \Phi (z) $ is continuous for $ z \in F $.
\end{lem}

\begin{pf}
	Given any two points $ z^{(1)}, z^{(2)} \in F $, and they are close to each other, i.e. $ \| z^{(1)} - z^{(2)} \| \leq \delta $. Let
	\begin{align*}
	& y^{(1)} = \underset{y}{argmin}\{P(y)\ |\ y\in F_{z^{(1)}} \} \; \text{and} \\
	& y^{(2)} = \underset{y}{argmin}\{P(y)\ |\ y\in F_{z^{(2)}} \},
	\end{align*}
	then we have $ \Phi (z^{(1)})=P(y^{(1)}) $ and $ \Phi (z^{(2)})=P(y^{(2)}) $. \\
	Without loss of generality, let $ \Phi (z^{(1)}) \leq \Phi (z^{(2)}) $, i.e. $ P(y^{(1)}) \leq P(y^{(2)}) $. From Lemma~\ref{lem:continuity}, the continuity of $ F_{z} $, there exists $ \hat{y}^{(2)} \in F_{z^{(2)}} $, such that $ \| \hat{y}^{(2)} - y^{(1)} \| \leq \eta $ for any $ \eta >0 $. Then since $ P(y) $ is continuous, we have
	\begin{equation} \label{eq:P_conti}
		| P(\hat{y}^{(2)}) - P(y^{(1)}) | \leq \epsilon.
	\end{equation}
	Since $ y^{(2)} $ is the minimizer of $ P(y) $ in $ F_{z^{(2)}} $, $ P(y^{(2)}) \leq P(\hat{y}^{(2)}) $. Therefore, $ P(y^{(1)}) \leq P(\hat{y}^{(2)}) $ by assumption. Now~\eqref{eq:P_conti} becomes $ P(\hat{y}^{(2)}) - P(y^{(1)}) \leq \epsilon $. Again, because $ P(y^{(2)}) \leq P(\hat{y}^{(2)}) $, we have $ P(y^{(2)}) - P(y^{(1)}) \leq \epsilon $, i.e. $ \Phi (z^{(2)}) - \Phi (z^{(1)}) \leq \epsilon $, which means $ \Phi (z) $ is continuous. \qed
\end{pf}

With the continuity of $ \Phi(z) $ in hand, now we are ready to present our final convergence results:

\begin{thm}[\textbf{Global convergence}]
	The sequence $ \{z^{(k)} \} $ generated by the successive procedure~\eqref{eq:sub_prob}	is in $ F $, and has limit point $ z^{*} $, at which the corresponding sequence $ \{P(z^{(k)})\} $ attains its minimum, $ P(z^{*}) $. \\
	More importantly, $ P(z^{*}) = \Phi(z^{*}) $, i.e. $ z^{*} $ is a local optimum of the penalty problem in~\eqref{eq:final_prob} convexified at $ z^{*} $.
\end{thm}

\begin{pf}
	From Lemma~\ref{lem:invariance}, we have $ z^{(k)} \in F_{z^{(k)}} \subseteq F $, then
	\begin{equation*}
	P(z^{(k+1)}) = \underset{y}{min}\{P(y)\ |\ y\in F_{z^{(k)}} \} \leq P(z^{(k)}),
	\end{equation*}
	because $ z^{(k)} $ is a feasible point to this convex optimization problem, while $ z^{(k+1)} $ is the optimum. Therefore, the sequence $ \{P(z^{(k)})\} $ is monotonically decreasing. \\
	Also since $ F_{z^{(k)}} \subseteq F $ for all $ k $, we have
	\begin{equation*}
	P(z^{(k+1)}) \geq \underset{y}{min}\{P(y)\ |\ y\in F \},
	\end{equation*}
	which means the sequence $ \{P(z^{(k)})\} $ is bounded from below. Then by the monotone convergence theorem, see e.g. \cite{rudin1964principles}, $ \{P(z^{(k)})\} $ converges to its infimum. Due to the compactness of $ F_{z^{(k)}} $, this infimum is attained by all the convergent subsequences of $ \{z^{(k)}\} $. Let $ z^{*} \in F $ be one of the limit points, then $ \{P(z^{(k)})\} $ attains its minimum at $ P(z^{*}) $, i.e.
	\begin{equation} \label{eq:infimum}
	P(z^{*}) \leq P(z^{(k)}), \quad \forall \; k = 0, 1,...
	\end{equation}
	
	To show $ P(z^{*}) = \Phi(z^{*}) $, we note that since $ z^{*} \in F $, $ z^{*} \in F_{z^{*}} $ by Lemma~\ref{lem:invariance}. Thus we have
	\begin{equation*}
	\Phi(z^{*}) \coloneqq \underset{y}{min}\{P(y)\ |\ y\in F_{z^{*}} \} \leq P(z^{*}).
	\end{equation*}
	Now for the sake of contradiction, we suppose $ \Phi(z^{*}) < P(z^{*}) $. Then by Lemma~\ref{lem:continuity_phi}, the continuity of $ \Phi(z) $, there exists a sufficiently large $ k $ such that $ \Phi(z^{k}) < P(z^{*}) $. Then we have $ P(z^{(k+1)}) = \Phi(z^{k}) < P(z^{*}) $, which contradicts~\eqref{eq:infimum}. Therefore $ P(z^{*}) = \Phi(z^{*}) $ holds. \qed
\end{pf}

%

\section{Numerical Results}
In this section, we present numerical results that apply the \texttt{SCvx} algorithm to an aerospace problem. Consider a multi--rotor vehicle with state at time $t_i$ given by $x_i = \left\lbrack p_i^T,v_i^T \right\rbrack^T$, where $p_i\in\mathbb{R}^3$ and $v_i\in\mathbb{R}^3$ represent vehicle position and velocity at time $t_i$ respectively. We assume that vehicle motion is adequately modeled by double integrator dynamics with constant time step $\Delta t$ such that $x_{i+1} = Ax_i + B(u_i+g)$, where $u_i\in\mathbb{R}^3$ is the control at time $t_i$, $g\in\mathbb{R}^3$ is a constant gravity vector, $A$ is the discrete state transition matrix, and $B$ utilizes zero order hold integration of the control input.
Further, we impose a speed upper--bound at each time step, $\|v_i\|_2\leq V_{\max}$, an acceleration upper--bound such that $\|u_i\|_2\leq u_{\max}$ (driven by a thrust upper--bound), and a thrust cone constraint $\hat{n}^Tu_i \geq \|u_i\|_2\cos(\theta_{cone})$ that constrains the thrust vector to a cone pointing towards the unit vector $\hat{n}$ (pointing towards the ceiling) with angle $\theta_{cone}$. Finally, the multi--rotor must avoid a known set of cylindrical obstacles with the $j^{th}$ cylinder having center $p_{c,j}\in\mathbb{R}^2$ and radius $r_j\in\mathbb{R}_+$. Therefore, each state must satisfy the non--convex constraint given by $\|Hx_i-p_{c,j}\|_2\geq r_j$ where $H$ is a linear mapping that maps $p_i$ to its projection on the ground plane.

Given these constraints, the objective is to find a minimum fuel trajectory from a prescribed $x(t_0)$ to a known $x(t_f)$ with fixed final time $t_f$ and $N$ discrete points along with $n_{cyl}$ cylindrical obstacles to avoid:
\begin{equation}
\begin{aligned}
\operatorname{min}& \sum_{i=1}^N \|u_i\|\\
\textnormal{s.t. }& x_{i+1} = Ax_i + B(u_i+g),\ i=1,\ldots,N-1\\
& \|u_i\|_2 \leq u_{\max},\ i = 1,\ldots,N,\\
& \|v_i\|_2\leq V_{\max}, \ i = 1,\ldots,N,\\
& \hat{n}^Tu_i \geq \|u_i\|_2\cos(\theta_{cone}), \ i = 1,\ldots,N,\\
& \|Hx_i-p_{c,j}\|_2\geq r_j,\ i=1,\ldots,N,\ j=1,\ldots,n_{cyl},\\
& x_0 = x(t_0),\ x_N = x(t_f).
\end{aligned}
\end{equation}

\begin{table}
	\centering
	\caption{Parameter Values}
	\label{tab:pars}
	
	\begin{tabular}{lc|lc}
		\hline\hline
		Par. & Value &  Par. & Value\rule{0pt}{2.6ex} \\
		\hline
		$N$ & 25 & $\epsilon$ & $1\times 10^{-6}$\rule{0pt}{2.6ex}\\
		$V_{\max}$ & 2 m/s & $u_{\max}$ & 13.33 m/s$^2$\\
		$g$ & $\lbrack 0,\,0,\,-9.81\rbrack^T$ m/s$^2$ & $\theta_{cone}$ & 30 deg \\
		$p_0$ & $\lbrack -8,\,-1,\,0\rbrack^T$ m & $p_f$ & $\lbrack 8,\,1,\,0.5\rbrack^T$ m\\
		$v_0$ & $\lbrack 0,\,0,\,0\rbrack^T$ m/s & $v_f$ & $\lbrack 0,\,0,\,0\rbrack^T$ m/s\\
		$p_{c,1}$ & $\lbrack -1,\,0\rbrack^T$ m & $p_{c,2}$ & $\lbrack 4,\,-1 \rbrack^T$ m\\
		$r_1$ & 3 m & $r_2$ & 1.5 m\\
		$\lambda$ & 0 & $\hat{n}$ & $\lbrack 0,\,0,\,1\rbrack^T$\\
		$t_f$ & 15 s &&\\
		\hline\hline
	\end{tabular}
\end{table}

\begin{figure}[!ht]
	\centering
	\includegraphics[width=.5\linewidth,trim=20 0 30 15,clip]{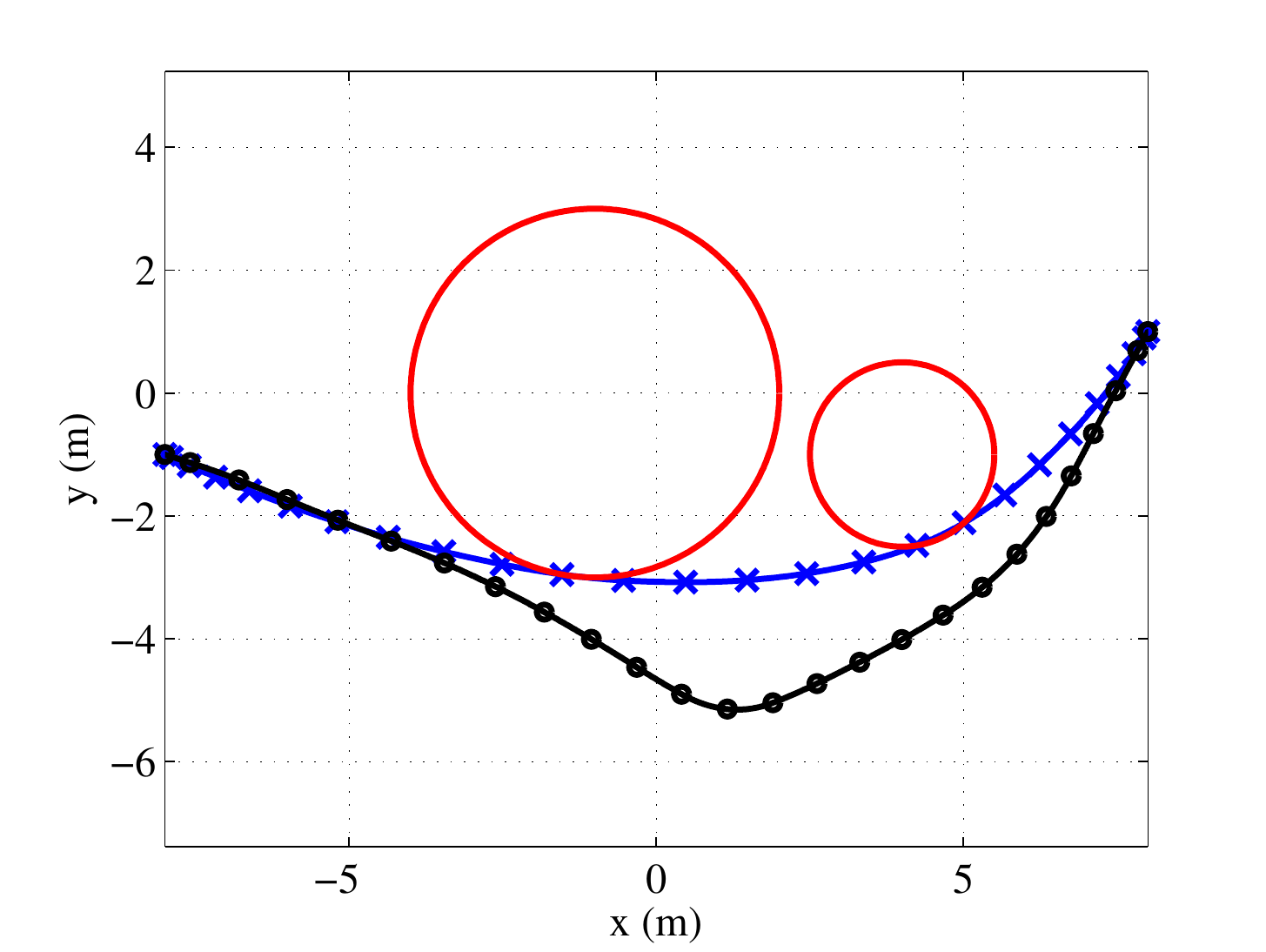}\label{fig:feas}
    \vspace{-2mm}
	\caption{2-D Ground plane projection of initial trajectory (black circles) and the converged trajectory (blue x's).\vspace{-1mm}}
	\label{fig:gps}
\end{figure}

The parameters given in Table \ref{tab:pars} are used to obtain the numerical results presented herein. An initial feasible trajectory is obtained by using the Trust Region Method (TRM) given in \cite{SCvx_cdc16} to solve the feasibility problem associated with this optimization. The ground plane projection of the feasible initial trajectory, $z^{\{0\}}$ is shown in Figure \ref{fig:feas} (black circles). Note that $z^{\{0\}}$ is not only feasible for the non--convex state constraints, but also for the convex constraints imposed on the trajectory.

\begin{figure}[!ht]
	\begin{center}
		\includegraphics[width=0.25\textwidth,trim=65 0 85 30,clip]{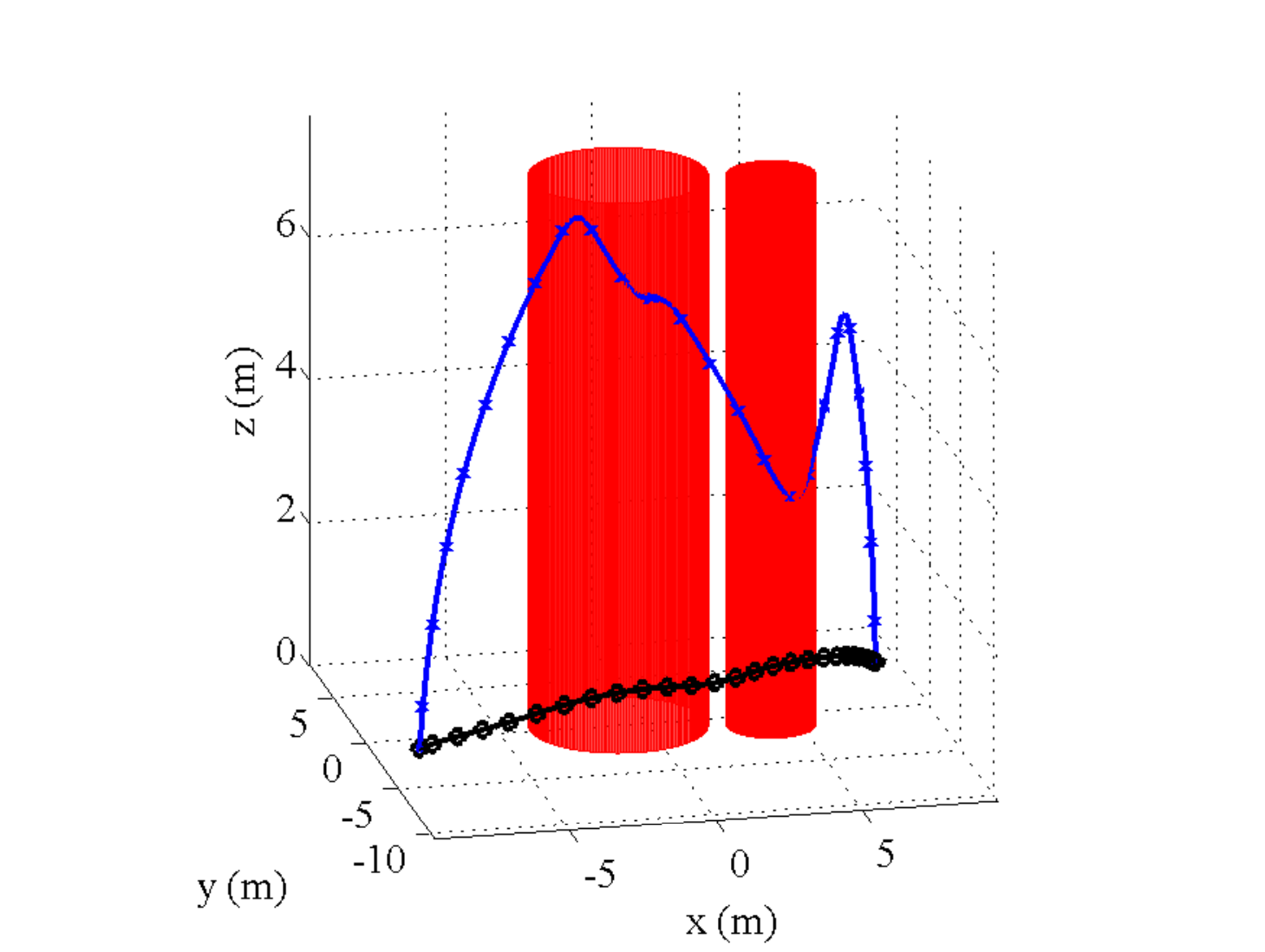}
        \vspace{-2mm}
		\caption{3--D plot of the initial trajectory (black circles) and the converged trajectory (blue x's). \vspace{-4mm}} 
		\label{fig:3d}
	\end{center}
\end{figure}

The \texttt{SCvx} algorithm is initiated with $z^{\{0\}}$, and is considered to have converged when the improvement in the cost of the linearized problem is less than $\epsilon$. Figure \ref{fig:3d} illustrates the converged trajectory that avoids the cylindrical obstacles while satisfying its actuator and mission constraints (blue x's). Note that the ground plane projection of the converged trajectory $z^{\{5\}}$ in Figure \ref{fig:gps} (blue x's) is different than that of $z^{\{0\}}$ and is characterized by having a smooth curve. For $z^{\{0\}}$, $\sum_{i=1}^N \|u_i\| = 269.01$ and at $z^{\{5\}}$, we have that $\sum_{i=1}^N \|u_i\| = 245.38$, so the cost of the converged trajectory is lower than that of the initial trajectory. At each iteration, the \texttt{SCvx} algorithm solves an SOCP, and therefore 5 SOCPs were solved in order to produce these results (in addition to 2 SOCPs for finding a feasible starting point). 

For comparison, we solve the same problem using only the TRM given in \cite{SCvx_cdc16}. The converged trajectories for both algorithms are identical to within $\epsilon$. The number of iterations ($k_{TRM}$ and $k_{\texttt{SCvx}}$) and time elapsed using each method ($t_{TRM}$ and $t_{\texttt{SCvx}}$) are reported in Table \ref{tab:times}. The two iterations necessary to find a feasible starting point for the \texttt{SCvx} algorithm are also reported. All times were found using the ECOS solver (\cite{alexd}) on a standard workstation with an Intel Xeon processor at 3.40 GHz and 16 GB of RAM.

\begin{table}[!ht]
	\centering
	\caption{Runtimes and Iterations}
	\label{tab:times}
	\setlength{\tabcolsep}{4pt}
	\begin{tabular*}{0.45\textwidth}{l|cccc|cc}
		\hline\hline
		Method & $k_{TRM}$ & $t_{TRM}$ & $k_{\texttt{SCvx}}$ & $t_{\texttt{SCvx}}$ & $k_{tot}$ & $t_{tot}$\\
		\hline
		TRM & 14 & 149.9 ms & - & - & 14 & 149.9 ms\\
		\texttt{SCvx} & 2 & 39.3 ms & 5 & 26.5 ms & 7 & 65.7 ms\\
		\hline\hline
	\end{tabular*}
\end{table}

\section{Conclusion}

The proposed successive convexification (\texttt{SCvx}) algorithm with its \textit{project-and-linearize} procedure solves a class of non-convex optimal control problems by solving a sequence of convex optimization problems. Further, we give a convergence proof demonstrating that the algorithm always converges to a local optimum. Numerical results suggest that this convergence process takes far fewer iterations than most trust-region-method based algorithms. Finally, the runtimes presented here show that the method has the potential to be used in real--time applications.

Future work includes finding a lifting procedure such that the algorithm can use any initial guess and using subgradient theories to drop the differentiability assumptions -- thus enabling the use of polyhedra shaped obstacles.

\begin{ack}
This research was supported in part by the Office of Naval Research Grant No. N00014-16-1-2318 and by the National Science Foundation Grants No. CMMI-1613235 and CNS-1619729.
\end{ack}

\bibliography{SCvx_ifac}

\end{document}